\documentclass[reprint,aip,graphicx,cha,showpacs]{revtex4-1}

\usepackage{float,MnSymbol}

\usepackage{amsmath,amsfonts,bm,bbm}
\usepackage{epsfig}

\def \R{\mathbb R}

\def \P{\mathbb P}
\def \E{\mathbb E}

\usepackage{color}

\usepackage{multirow}

\newcommand{\e}{{\mathrm e}}

\newcommand{\x}{{\bf x}}

\newcommand{\M}{\mathcal{M}}

\renewcommand{\theequation}{\arabic{section}.\arabic{equation}}
\newcommand{\james}[1]{{\color{black}#1}}

\begin{document}

\title[Synchronization of stochastic hybrid oscillators]{Synchronization of stochastic hybrid oscillators driven by a common switching environment}
\author{Paul C. Bressloff and James MacLaurin}
\address{Department of Mathematics, University of Utah, Salt Lake City, UT 84112 USA}

\date{\today}

\begin{abstract}
Many systems in biology, physics and chemistry can be modeled through ordinary differential equations, which are piecewise smooth, but switch between different states according to a Markov jump process. In the fast switching
limit, the dynamics converges to a deterministic ODE. In this paper we suppose that this limit ODE supports a stable limit cycle.
We demonstrate that a set of such oscillators can synchronize when they are uncoupled, but they share the same switching Markov jump process. The latter is taken to represent the effect of a common randomly switching environment. We determine the leading order of the Lyapunov coefficient governing the rate of decay of the phase difference in the fast switching limit. The analysis bears some similarities to the classical analysis of synchronization of stochastic oscillators subject to common white noise. However the discrete nature of the Markov jump process raises some difficulties: in fact we find that the Lyapunov coefficient from the quasi-steady-state approximation differs from the Lyapunov coefficient one obtains from a second order perturbation expansion in the waiting time between jumps. Finally, we demonstrate synchronization numerically in the radial isochron clock model and show that the latter Lyapinov exponent is more accurate.
\end{abstract}
\maketitle 

\noindent {\bf There are a growing number of systems in physics and biology where a population of oscillators can be synchronized by a randomly fluctuating external input applied globally to all of the oscillators, even if there are no interactions between the oscillators (noise-induced phase synchronization). Experimental evidence for such an effect has been found in neural oscillations of the olfactory bulb, synthetic genetic oscillators, laser dynamics, and variations in geographically separated animal populations. Most previous studies of noise-induced phase synchronization have taken the oscillators to be driven by common Gaussian noise. Typically, phase synchronization is established by constructing the Lyapunov exponent for the dynamics of the phase difference between a pair of oscillators and averaging with respect to the noise. If the averaged Lyapunov exponent is negative definite, then the phase difference is expected to decay to zero in the large time limit, establishing phase synchronization. In this paper we extend the theory of noise-induced synchronization to the case of a common randomly switching environment. Each oscillator then evolves according to a piecewise deterministic Markov process, which involves the coupling between a piecewise continuous differential equation and a time-homogeneous Markov chain. In the fast switching
limit, the dynamics converges to a deterministic ODE, which is assumed to support a stable limit cycle. 
We demonstrate that an uncoupled population of such oscillators can synchronize when they share the same switching Markov jump process. We determine the leading order of the Lyapunov coefficient governing the rate of decay of the phase differences in the weak noise regime (fast but finite switching rates), and show that it differs from the standard expression obtained using a Gaussian approximation of the noise.

}

\rule{6cm}{0.4pt}

\section{Introduction}

Self--sustained oscillations in biological, physical and chemical systems are often described in terms of limit cycle oscillators where the timing along each limit cycle is specified in terms of a single phase variable. Phase reduction methods can then be used to analyze synchronization of an ensemble of weakly-coupled oscillators by approximating the high--dimensional limit cycle dynamics as a closed system of equations for the corresponding phase variables \cite{Winfree80,Kuramoto84,Erm84,Glass88,Erm91,Erm96,Holmes04,Ashwin16,Nakao16}. More recently, there has been considerable interest in applying phase reduction methods to the analysis of noise-induced phase synchronization \cite{Teramae04,Goldobin05,Nakao07,Yoshimura08,Teramae09,Ly09,Erm10}. This concerns the observation that a population of oscillators can be synchronized by a randomly fluctuating external input applied globally to all of the oscillators, even if there are no interactions between the oscillators. Evidence for such an effect has been found in experimental studies of neural oscillations in the olfactory bulb \cite{Galan08}, and the synchronization of synthetic genetic oscillators \cite{Zhou08,Paulsson16}. A related phenomenon is the reproducibility of a dynamical system's response when repetitively driven by the same fluctuating input, even though initial conditions vary across trials. Examples include the spike-time reliability of single neurons  \cite{Mainen95,Galan06}, improvements in the reproducibility of laser dynamics \cite{Uchida04}, and synchronized variations in wild animal populations located in distinct, well-separated areas caused by common environmental 
fluctuations \cite{Grenfell98}. 

\begin{figure}[b!]
\begin{center}
\includegraphics[width=8cm]{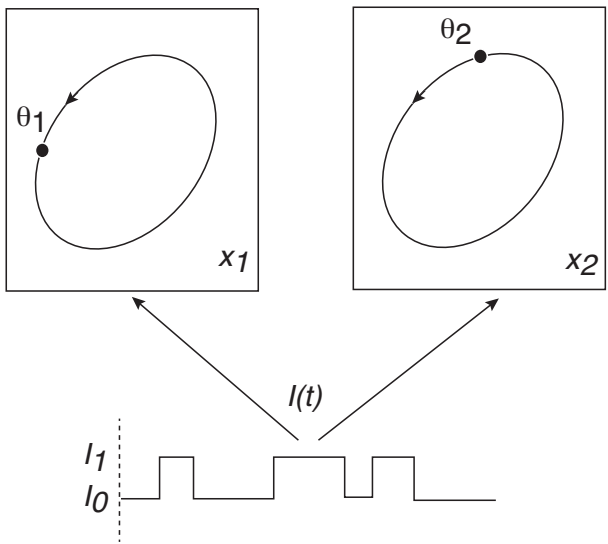}
\caption{\small Pair of non-interacting limit cycle oscillators with phases $\theta_j(t)$, $j=1,2$, driven by a common switching external input $I(t)$}
\label{dich}
\end{center}
\end{figure}

Most studies of noise-induced synchronization take the oscillators to be driven 
 by common Gaussian noise. Typically, phase synchronization is established by constructing the Lyapunov exponent for the dynamics of the phase difference between a pair of oscillators and averaging with respect to the noise. If the averaged Lyapunov exponent is negative definite, then the phase difference is expected to decay to zero in the large time limit, establishing phase synchronization.
However, it has also been shown that common Poisson-distributed random impulses, dichotomous or telegrapher noise, and other types of noise generally induce synchronization of limit-cycle oscillators \cite{Nakao05,Nagai05,Goldobin10}. Consider, in particular, the case of an additive dichotomous noise signal $I(t)$ driving a population of $M$ identical non-interacting oscillators according to the system of equations $\dot{x}_j=F(x_j)+I(t)$, where $x_j\in \R^d$ is the state of the $j$th oscillator, $j=1,\ldots,M$ \cite{Nagai05}, see Fig. \ref{dich}. Here $I(t)$ switches between two values $I_0$ and $I_1$ at random times generated by a two-state Markov chain \cite{Bena06}. That is, $I(t)=I_0(1-N(t))+I_1N(t)$ for $N(t)\in \{0,1\}$, with the time $T$ between switching events taken to be exponentially distributed with mean switching time $\tau$. Suppose that each oscillator supports a stable limit cycle for each of the two input values $I_0$ and $I_1$. It follows that the internal state of each oscillator randomly jumps between the two limit cycles. Nagai et al \cite{Nagai05} show that in the slow switching limit (large $\tau$), the dynamics can be described by random phase maps. Moreover, if the phase maps are monotonic, then
the associated Lyapunov exponent is generally negative and phase synchronization is stable.

The dichotomous noise-driven system is just one example of a much more general class of randomly switching processes known as {\em piecewise deterministic Markov processes} (PDMPs) \cite{Davis84,Bressloff17}. More explicitly, let $N(t) \in \Gamma \equiv \{0,\cdots,\Lambda_0-1\}$ denote the state of the randomly switching environment. When the environmental state is $N(t)=n$, each oscillator $x_i(t)$ evolves according to the piecewise deterministic ordinary differential equation (ODE)
$\dot{x_i}=F_n(x_i),\quad i=1,\ldots,M$,
where the vector  field $F_n: \R^d \to \R^d$ is a smooth function. The discrete stochastic variable $N(t)$ evolves according to a stationary, continuous-time Markov chain with transition matrix ${\bf W}$. The additive dichotomous noise case is recovered by taking $\Lambda_0=2$ and $F_n(x)=F(x)+I_n$. 
One major application of PDMPs is to stochastic gene regulatory networks, where the continuous variables $x_j$ are the concentrations of protein products (and possibly mRNAs) and the discrete variables represent the various activation/inactivation states of the genes \cite{Kepler01,Bose04,Zeiser08,Smiley10,Paulsson11,Newby12,Singh14,Koeppl14,Newby15,Hufton16,Levien18}. The common randomly switching environment could represent the state of a promoter site that is common to a pair of genes within the same cell, or the state of the extracellular environment that drives gene expression in a population of cells. \james{It} is thought that synchronous oscillations may have an important functional purpose in systems biology \cite{Paulsson11}

In this paper we develop the theory of noise-induced synchronization for a population of non-interacting PDMPs evolving under a common randomly switching environment. (The population model is presented in section II). In the slow switching limit one could generalize the approach of Nagai et al \cite{Nagai05} by assuming that each of the vector fields $F_n(x_i)$, $n\in \Gamma$,  supports a stable limit cycle and constructing the associated random phase maps. Here, instead, we consider the fast switching regime in which
 the transition rates between the discrete states $n\in \Gamma$ are much faster than the relaxation rates of the piecewise deterministic dynamics for $x_i\in \R^d$. Thus there is a separation of time scales between the discrete and continuous processes, so that if $R$ is the characteristic relaxation rate of the continuous dynamics, then $r/\epsilon$ is the characteristic transition rate of the Markov chain for some small positive dimensionless parameter $\epsilon$. If the Markov chain is ergodic, then in the fast switching or adiabatic limit $\epsilon\rightarrow 0$ one obtains a deterministic dynamical system in which one averages the piecewise dynamics with respect to the corresponding unique stationary distribution. Suppose that in the deterministic limit we have a population of independent limit cycle oscillators. Since there is no coupling or remaining external drive to the oscillators, their phases are uncorrelated. The basic issue we wish to address is whether or not phase synchronization occurs when $\epsilon >0$; we will refer to the resulting oscillators as {\em stochastic hybrid limit cycle oscillators}. We will proceed by constructing the Lyapunov exponent for a pair of such oscillators driven by a common randomly switching environment. 
 
In section III we obtain an approximate expression for the Lyapunov exponent by considering a quasi-steady-state (QSS) diffusion approximation of the underlying PDMPs \cite{Newby10a} (see also appendix A), in which each oscillator is approximated by a stochastic differential equation (SDE) with a common Gaussian input. This allows us to adapt previous work on the phase reduction of stochastic limit cycle oscillators \cite{Teramae04,Goldobin05,Nakao07,Teramae09}, and thus establish that phase synchronization occurs under the diffusion approximation. However, the QSS approximation is only intended to be accurate over timescales that are longer than $O(\epsilon)$. Hence, it is unclear whether or not the associated Lyapunov exponent is accurate, since it is obtained from averaging the fluctuations in the noise over infinitesimally small timescales. Therefore, in section IV we derive a more accurate expression for the Lyapunov exponent by working directly with an exact implicit equation for the phase dynamics. We exploit the fact that multiple switching events (jumps) occur during small excursions around the limit cycle for small $\epsilon$, which allows us to express the Lyapunov exponent in terms of discrete sums over these events. Taking expectations then yields an expression for the Lyapunov exponent that differs significantly from the one obtained using the diffusion approximation. Note, however, that both are negative definite, so they both imply phase synchronization but at different rates. Our derivation of the Lyapunov exponent from the exact phase equation also allows us to obtain greater insights into the nature of the QSS approximation and the meaning of the associated Brownian motion (see also appendix B). Finally, we illustrate the theory by considering the particular example of radial isochron clocks (section V).

\section{Population of stochastic hybrid limit cycle oscillators}

Consider a population of identical, noninteracting dynamical systems labeled $i=1,\ldots,M$, whose states are described by the pair
$(x_i(t),N(t)) \in \Sigma \times\Gamma$, where $x_i(t)$ is a continuous variable in a connected bounded domain $\Sigma\subset \R^d$ and $N(t)$ is an $i$-independent discrete stochastic variable taking values in the finite set $\Gamma \equiv \{0,\cdots,\Lambda_0-1\}$. The latter represents the state of an environment that is common to all members of the population.
When the environmental state is $N(t)=n$, $x_i(t)$ evolves according to the piecewise deterministic ODE
\begin{equation}
\label{pdmp}
\dot{x_i}=F_n(x_i),\quad i=1,\ldots,M,
\end{equation}
where the vector  field $F_n: \R^d \to \R^d$ is a smooth function. We assume that the dynamics of $x_i$ is confined to the domain $\Sigma$. The discrete stochastic variable is taken to evolve according to a homogeneous, continuous-time Markov chain with $x$-independent generator ${\bf A}$,  where
\[A_{nm}=W_{nm}-\delta_{n,m}\sum_{k\in \Gamma} W_{kn},\]
and ${\bf W}$ is the transition matrix. We make the further assumption that the chain is irreducible, that is, there is a non-zero probability of transitioning, possibly in more than one step, from any state to any other state of the Markov chain. This implies the existence of a unique invariant probability distribution on $\Gamma$, denoted by $\rho$, such that
\begin{equation}
\label{Wstar}
\sum_{m\in \Gamma}A_{nm}\rho_m=0,\quad \forall n \in \Gamma.
\end{equation}
As a simple example, suppose that $N(t)$ evolves according to a two-state Markov chain. That is, $N(t)\in \Gamma \equiv \{0,1\}$ and the generator of the Markov chain is given by the matrix
\begin{equation}
{\bf A}=\left (\begin{array}{cc} -k_- &k_+  \\ k_- & -k_+  \end{array} \right ).
\label{A}
\end{equation}
The corresponding stationary distribution of the Markov chain then has components
\begin{equation}
\rho_0 =\frac{k_+ }{k_++k_- },\quad \rho_1= \frac{k_- }{k_+ +k_- }.
\end{equation}

\begin{figure}[b!]
\begin{center}
\includegraphics[width = 8cm]{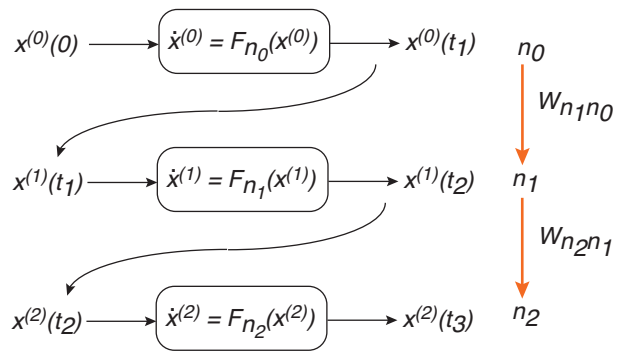}
\caption{Schematic diagram of a PDMP for a sequence of jump times $\{t_1,\ldots\}$, and a corresponding set of discrete states $\{n_0,n_1,\ldots\}$. See text for details}
\label{fig1}
\end{center}
\end{figure}

Eq. (\ref{pdmp}) defines a PDMP \cite{Davis84} on $\R^{d}$ for each $i=1,\ldots,M$, also known as a stochastic hybrid system (SHS). A useful way to implement the PDMP is to decompose the transition matrix of the Markov chain as $W_{nm}=P_{nm}\lambda_m$,
with $\sum_{n\neq m}P_{nm}=1$.  Here $\lambda_m$ is the rate of the exponential waiting time density for transitions from state $m$, whereas $P_{nm}$ is the probability of the transition $m\rightarrow n$, $n\neq m$. Suppose that $N(t)=n_0$ and let $t_1$ be an exponentially distributed random variable with rate $\lambda_{n_0}$. That is, 
\[\P(t_1 < t) =1- \exp \left( -\lambda_{n_0}t\right).\]
Then in the random time interval $s\in [0,\,t_1)$ the state of the $i$th system is $(x^{(0}(s),n_0)$ with $x^{(0)}$ evolving according to Eq. (\ref{pdmp}) for $n=n_0$. (For the moment, we drop the population label $i$.) At the random time $t_1$ we choose an internal state $n_1 \in \Gamma$ with probability $P_{n_1n_0}$, and call $x^{(1)}(t)$ the solution of the following Cauchy problem on $[t_1,\infty)$:
\[
\left\{
\begin{array}{lcl}
\dot{x}^{(1)}(t) & = & F_{n_1}(x^{(1)}(t),\quad t \geq t_1\\ 
x^{(1)}(t_1) & = & x^{(0)}(t_1)
\end{array}
\right.
\]
Iterating this procedure, we construct a sequence of increasing jumping times $(t_k)_{k \geq 0}$ (setting $t_0=0$) and a corresponding sequence of internal states $(n_k)_{k \geq 0}$. The evolution $(x(t),\,N(t))$ is then defined as, see Fig. \ref{fig1}
\begin{equation}
\label{hs}
(x(t),N(t))=(x^{(k)}(t),n_k) \quad \mbox{if}\ t_k \leq t <t_{k+1}.
\end{equation}

 Introduce the population vector ${\bf x}(t)=(x_1(t),\ldots,x_M(t))$ and define the probability density $p_n(\x,t) $, given the initial conditions $x_i(0)=x_{i,0},N(0)=n_0$, according to
\[p_n(\x,t)d\x=\P\{\x(t)\in (\x,\x+d\x),\, N(t) =n|\x_0,n_0\}.\]
It can be shown that $p$ evolves according to the forward differential Chapman-Kolmogorov (CK) equation \cite{Gardiner09,Bressloff17}
\begin{equation}
\label{CKH}
\frac{\partial p_n}{\partial t}={\mathbb L}p_n ,
\end{equation}
with the generator ${\mathbb L}$ defined according to
\begin{equation}
\label{linH}
{\mathbb L} p_n(\x,t)=-\sum_{i=1}^M \nabla_i \cdot \left [F_n(x_i)p_n(\x,t)\right ] +\frac{1}{\epsilon}\sum_{m\in \Gamma}A_{nm}p_m(\x,t).
\end{equation}
Here $\nabla_i$ denotes the $d$-dimensional gradient operator with respect to $x_i$.
The first term on the right-hand side represents the probability flow associated with the piecewise deterministic dynamics for a given $n$, whereas the second term represents jumps in the discrete state $n$. Note that we have rescaled the matrix ${\bf A}$ by introducing the dimensionless parameter $\epsilon$, $\epsilon >0$. This is motivated by the observation that many applications of PDMPs involve a separation of time-scales between the relaxation time for the dynamics of the continuous variables $x$ and the rate of switching between the different discrete states $n$ of the environment. The fast switching limit then corresponds to the case $\epsilon \rightarrow 0$. Now introduce the averaged vector field $\overline{F}: \R^d \to \R^d$ by
\begin{equation}
\label{Fav}
\overline{F}(x)=\sum_{n \in \Gamma}\rho_n F_n(x)
\end{equation}
and define the averaged system
\begin{equation}
\label{mft}
\left\{
\begin{array}{lcl}
\dot{x_i}(t) & = &  \overline{F}(x_i(t)),\quad i=1,\ldots,M,\\
x_i(0) & = & x_0
\end{array}
\right.
\end{equation}
Intuitively speaking, one expects the PDMP (\ref{pdmp}) to reduce to the deterministic dynamical system (\ref{mft}) in the fast switching limit $\epsilon \rightarrow 0$. That is, for sufficiently small $\epsilon$, the Markov chain undergoes many jumps over a small time interval $\Delta t$ during which $\Delta x\approx 0$, and thus the relative frequency of each discrete state $n$ is approximately $\rho_n$. This can be made precise in terms of a law of large numbers for PDMPs proven in \cite{Faggionato10}. 

In the fast switching (deterministic) limit, each member of the population becomes independent, since the dependence on the current state of the environment disappears. In this paper, we will assume that for each $i=1,\ldots,M$, the averaged dynamical system (\ref{mft}) supports a set of stable periodic solutions with the same natural frequency $\overline{\omega}=2\pi/\overline{\Delta}$. That is, we have a population of identical, independent oscillators in the fast switching limit. In state space, each periodic solution is an isolated attractive trajectory
or limit cycle.
The dynamics on the limit cycle can be described by a uniformly rotating phase such that
\begin{equation}
\frac{d\theta_i}{dt}=\overline{\omega},
\end{equation}
and $x_i={\Phi}(\theta_i(t))=\Phi(\overline{\omega} t+\psi_i)$ with ${\Phi}$ a $2\pi$-periodic function and $\psi_i$ the initial phase. Note that $\Phi$ satisfies the equation
\begin{equation}
\label{fi}
\overline{\omega} \frac{d\Phi}{d\theta} = \overline{F}(\Phi(\theta)).\end{equation}
Differentiating both sides with respect to $\theta$ gives
\begin{equation}
\frac{d}{d\theta} \left (\frac{d\Phi}{d\theta}\right )=\overline{\omega}^{-1}\overline{J}( \theta)\cdot \frac{d\Phi}{d\theta},
\label{nadj}
\end{equation}
where
$\overline{J}$ is the $2\pi$-periodic Jacobian matrix 
\begin{equation}
\label{Jac}
\overline{J}_{ab}(\theta)\equiv \left . \frac{\partial \overline{F}_a}{\partial x_b}\right |_{x=\Phi(\theta)}
\end{equation}
for $a,b=1,\ldots,d$.

In the deterministic limit, there is no mechanism for phase synchronizing the population of oscillators, since $\theta_j(t)-\theta_i(t)=\psi_j-\psi_i$ for all $t$. The main issue we wish to address in this paper is whether or not the presence of a common switching environment can synchronize the population of stochastic hybrid oscillators when $\epsilon >0$, analogous to the noise-driven synchronization of SDEs \cite{Teramae04,Goldobin05,Nakao07,Yoshimura08,Teramae09}.

\setcounter{equation}{0}
\section{Diffusion approximation and phase SDE} 

One approach to analyzing synchronization in the fast switching regime ($0<\epsilon \ll 1$) is to use a QSS diffusion or adiabatic approximation, in which the CK Eq. (\ref{CKH}) is approximated by a Fokker-Planck (FP) equation for the total density $C(\x,t)=\sum_n p_n(\x,t)$ \cite{Newby10a}. The latter determines the probability distribution of solutions of a corresponding SDE for a population of oscillators driven by a common $O(\sqrt{\epsilon})$ multiplicative noise term, which can then be reduced to an effective SDE for the phases along the lines of Ref. \cite{Teramae04}. The resulting FP equation in the Stratonovich representation takes the form (see appendix A)
\begin{subequations}
\begin{align}
\label{zFP0}
\frac{\partial C}{\partial t}&=- \sum_{i=1}^M \nabla_i \cdot \left [\overline{F}(x_i) C\right ]\\
&-\epsilon \sum_{i,j=1}^M\sum_{m,n\in \Gamma}A_{mn}^{\dagger}\rho_n \nabla_i\cdot  \left [ G_{m}(x_i) \nabla_j \cdot (G_n(x_j)C ) \right ],\nonumber
\end{align}
with 
\begin{eqnarray}
\label{G}
G_{m}(x) &=&F_{m}(x) -\overline{F}(x) .
\end{eqnarray}
\end{subequations}
Since Eq. (\ref{zFP0}) is symmetric with respect to the exchange $(ni,mj)\leftrightarrow (mj,ni)$ we can replace ${\bf A}\rho$ by its symmetric part 
\begin{equation}
\tilde{A}_{mn} =\frac{1}{2}(A^{\dagger}_{mn} \rho_n  + A^{\dagger}_{nm} \rho_m).\end{equation}
 Note that the matrix $\tilde{\bf A}$ is negative definite, as we demonstrate in the Appendix. For example, in the case of a two-state Markov chain, $\tilde{\bf A}=\mbox{diag}(-\rho_0/(k_++k_-),-\rho_1/(k_++k_-))$

It follows that under the diffusion approximation, the PDMP (\ref{pdmp}) can be approximated by the Stratonovich SDE
\begin{equation}
dX_i(t)=  \overline{F}(X_i)dt +  \sqrt{2\epsilon}\sum_{m,n\in \Gamma} G_m(X_i)B_{mn}dW_n(t)\label{dep2}
\end{equation}
for $i=1,\ldots,M$,
where ${\bf B}{\bf B}^{\top}=-\tilde{\bf A}$, and $W(t)$ is a vector of uncorrelated Brownian motions in $ \mathbb{R}^{M}$,
\[
\mathbb{E}\big[{\bf W}(t){\bf W}(t)^{\top} \big] = t {\bf I},
\]
and ${\bf I}$ is the identity matrix. (Note that it doesn't matter which Hermitian square root of $-\tilde{\bf A}$ we take for ${\bf B}$, since they all yield the same statistical behavior of $X_i(t)$.)
Eq. (\ref{dep2}) represents a population of independent, non-interacting limit cycle oscillators, driven by common external white noise. One can now use phase reduction methods developed for SDEs.

\subsection{Phase reduction}

First, suppose that
the noise amplitude $\epsilon$ is sufficiently small relative to the rate of attraction to the limit cycle, so that deviations transverse to the
limit cycle are also small (up to some exponentially large stopping time). This suggests that the definition of a phase variable persists in the stochastic setting, and one can derive a stochastic phase equation by decomposing the solution to the SDE (\ref{dep2}) according to
\begin{equation}
\label{deco}
X_i(t)=\Phi(\theta_i(t))+\sqrt{\epsilon}v_i(t),
\end{equation}
with $\theta_i(t)$ and $v_i(t)$ corresponding to the phase and amplitude components, respectively, of the $i$th oscillator. 
However, there is not a unique way to define the phase $\theta_i$, which reflects the fact that there are different ways of projecting the exact solution onto the limit cycle \cite{Gonze02,Koeppl11,Bonnin17,Maclaurin18a}, see Fig. \ref{isochrone}. One well-known approach is to use the method of isochrons \cite{Teramae04,Nakao07,Yoshimura08,Teramae09}, which we briefly outline here. 

\begin{figure}[b!]
\begin{center}
\includegraphics[width=7cm]{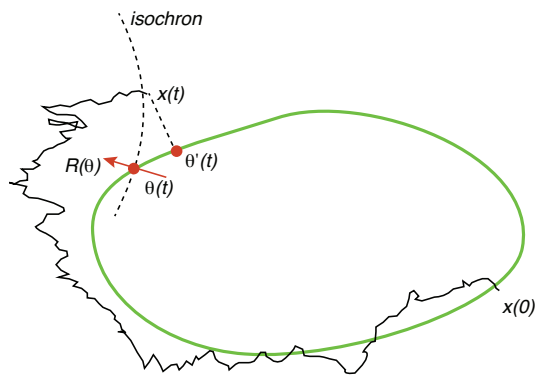}
\caption{\small Different choices of amplitude-phase decomposition. Two possibilities are orthogonal projection with phase $\theta'(t)$ and isochronal projection with phase $\theta(t)$. In the latter case, the response to perturbations depends on the phase response curve ${R}(\theta)$, which is normal to the isochron at the point of intersection with the limit cycle.}
\label{isochrone}
\end{center}
\end{figure}

Consider the unperturbed deterministic system $\dot{x}_i=\overline{F}(x_i)$. 
Stroboscopically observing the system at time intervals of length $\overline{\Delta}$ leads to a
Poincare mapping
\begin{equation*}
x_i(t)\rightarrow x_i(t+\overline{\Delta})\equiv {\mathcal P}(x_i(t)),
\end{equation*}
for which all points on the limit cycle are fixed points. Choose a point $x_i^*$ on the limit cycle and consider
all points in the vicinity of $x_i^*$ that are attracted to it under the action of ${\mathcal P}$. They form a
$(d-1)$-dimensional hypersurface ${\mathcal I}$ called an isochron  \cite{Winfree80,Kuramoto84,Glass88,Erm96,Holmes04}, crossing the limit cycle at $x_i^*$. A unique isochron can be drawn through each point on the limit cycle (at least locally) so the isochrons can be parameterized by the phase,
${\mathcal I}={\mathcal I}(\theta_i)$. Finally, the definition of phase is extended by taking all points $x_i\in {\mathcal I}(\theta_i)$ to have the same phase,
$\Theta(x_i)=\theta_i$, which then rotates at the natural frequency $\overline{\omega}$ (in the
unperturbed case). Hence, for an unperturbed oscillator in the vicinity of the limit cycle we have
\begin{eqnarray}
\overline{\omega} = \frac{d\Theta(x_i)}{dt}=\nabla \Theta(x_i)\cdot \frac{dx_{i}}{dt} =\nabla \Theta(x_i)\cdot F(x_i) .\nonumber
\end{eqnarray}

Now consider the Stratonovich SDE (\ref{dep2}). For the moment, we replace the $O(\epsilon^{1/2})$ term by a bounded deterministic function $H\in \R^d$ so that we have the perturbed deterministic equation
\[\frac{dx_i}{dt}=\overline{F}(x_i)+\sqrt{\epsilon}H(x_i,t).\]
The additional complications arising from a stochastic perturbation will be addressed below.
Differentiating the isochronal phase using the chain rule gives
\begin{align*}
\frac{d\Theta}{dt}&=\nabla \Theta(x_i)\cdot (\overline{F}(x_i)+\sqrt{\epsilon} H({x_i},t))\\
&=\overline{\omega}+\sqrt{\epsilon}\nabla \Theta(x_i)\cdot H(x_i,t).
\end{align*}
We now make the approximation that deviations of $x_i$ from the limit cycle are ignored on the right-hand side by setting $x_i(t)=\Phi(\theta_i(t))$ with $\Phi$ the $2\pi$-periodic solution on the limit cycle. This then yields  the closed phase equation
\begin{equation}
\frac{d\theta_i}{dt}=\overline{\omega}+\sqrt{\epsilon} \sum_{a=1}^d R_a(\theta_i)H_a(\Phi(\theta_i),t),
\label{cbphase}
\end{equation}
where
\begin{equation}
R_a(\theta)= \left . \frac{\partial \Theta}{\partial x_a}\right |_{x_a=\Phi(\theta)}
\label{Q2}
\end{equation}
is a $2\pi$-periodic function of $\theta$ known as the $a$th component of the {\em phase resetting curve} (PRC)
\cite{Winfree80,Kuramoto84,Glass88,Erm96,Holmes04}. One way to evaluate the PRC is to exploit the fact that it is the $2\pi$-periodic solution of the linear equation
\begin{equation}
\label{adj}
\overline{\omega}\frac{dR(\theta)}{d\theta}=-\overline{J}(\theta)^{\top}\cdot  R(\theta),
\end{equation}
under the normalization condition
\begin{equation}
\label{norry}
R(\theta)\cdot \frac{d\Phi(\theta)}{d\theta}=1.
\end{equation}
$\overline{J}(\theta)^{\top}$ is the transpose of the Jacobian matrix $\overline{J}(\theta)$.

Returning to the Stratonovich SDE (\ref{dep2}), treating the stochastic perturbation along identical lines to the deterministic case, (i.e. substituting $H(t)dt = \sqrt{2} \sum_{m,n\in \Gamma}\left [F_{m}(\Phi(\theta_i))-\overline{F}(\Phi(\theta_i))\right ]B_{mn}dW_n(t)$) and exploiting the fact that the usual rules of calculus hold (in contrast to Ito SDEs), would then lead to the following SDE for the phase $\theta_i(t)$:
\begin{align}
\label{dphase}
d\theta_i&=\overline{\omega} dt+\sqrt{2\epsilon} \sum_{m,n\in \Gamma}{\mathcal F}_m(\theta_i)B_{mn}dW_n(t),
\end{align}
where
\begin{equation}
\label{calF}
{\mathcal F}_n(\theta_i)=\sum_{a=1}^dR_a(\theta_i)\left [F_{n,a}(\Phi(\theta_i))-\overline{F}_a(\Phi(\theta_i))\right ].
\end{equation}
Introducing the population phase vector ${\bm \theta}(t)=(\theta_1(t),\ldots,\theta_M(t))$,
the corresponding phase FP equation for the population probability density $C({\bm \theta},t)$ is
\begin{align}
\label{FPphase}
\frac{\partial C}{\partial t}&=-\overline{\omega} \sum_{i=1}^M \frac{\partial C}{\partial \theta_i}\\
&+\epsilon \sum_{i,j=1}^M\sum_{m,n\in \Gamma}\tilde{A}_{mn}\frac{\partial}{\partial \theta_i} \left [ {\mathcal F}_{m}(\theta_i) \frac{\partial }{\partial \theta_j}({\mathcal F}_n(x_j)C ) \right ].\nonumber 
\end{align}
However, there are a number of major differences from the deterministic case. Probably the most significant is that a Wiener process is not bounded, which means that over sufficiently long time intervals there is a small but non-zero probability that the stochastic term induces large deviations from the limit cycle, resulting in a breakdown of the perturbation analysis. This issue can be addressed using variational methods and large deviation theory \cite{Giacomin16,Maclaurin18a}, which show that for sufficiently small $\epsilon$, the system remains in a neighborhood of the limit cycle up to exponentially long times. The second issue is that there are no $O(\sqrt{\epsilon})$ corrections to the deterministic part of the phase equation so that one has to go to $O(\epsilon)$ in order to determine the leading order corrections to the drift term. There are two sources of $O(\epsilon)$ terms: one arises from the coupling between the phase and amplitude fluctuations transverse to the limit cycle, and the second arises from changing between Stratonovich and Ito versions of the SDE based on Ito's formula \cite{Teramae09,Giacomin16,Bonnin17,Maclaurin18a}. The precise form of these terms will also depend on the particular choice of phase reduction method. However, if the limit cycle is sufficiently attracting then they tend to have a small effect \cite{Teramae09}. Moreover, such drift terms do not contribute to the leading order expression for Lyapunov exponent describing the evolution of phase differences, see below Eq. (\ref{Lyp1}). Therefore, we shall drop such contributions in our subsequent analysis, and reinterpret Eq. (\ref{dphase}) as an Ito SDE.  

\subsection{Lyapunov exponent}

We now use the phase SDE (\ref{dphase}) interpreted in the Ito sense to investigate the effects of a common switching environment in the small $\epsilon$ regime, following previous studies \cite{Teramae04,Nakao07,Yoshimura08,Teramae09,Ly09}. The first step is to consider the SDE for the phase difference $\psi=\theta_j-\theta_1$ for any fixed $j$ such that $j\neq 1$. Assuming $\psi$ is infinitesimally small, we have
\begin{equation}
\frac{d\psi }{dt}=\sqrt{2\epsilon} \psi(t)\left [\sum_{m,n\in \Gamma}{\mathcal F}'_m(\phi)B_{mn}dW_n(t)\right ],
\end{equation}
where $'$ denotes differentiation with respect to $\phi$ and we have set $\theta_1(t)=\phi(t)$ so that $\phi(t)$ evolves according to Eq. (\ref{dphase}) for $i=1$. Introducing a new variable $y=\log(\psi)$ and using Ito's formula yields the SDE
\begin{align}
dy&=-\epsilon \sum_{m,n\in \Gamma}{\mathcal F}'_m(\phi)\tilde{A}_{mn}{\mathcal F}'_n(\phi)dt\nonumber \\
&\quad +\sqrt{2\epsilon} \left [\sum_{m,n\in \Gamma}{\mathcal F}'_m(\phi)B_{mn}dW_n(t)\right ].\label{eq: QSS log phase}
\end{align}
Define the Lyapunov exponent according to
\begin{align}
\lambda_{\rm QSS} &=\lim_{T\rightarrow \infty} \frac{1}{T} (y(T)-y(0)).\end{align}
(More precisely, we only take the limit in $T$ up until the time that the system leaves a neighborhood of the limit cycle. In previous work \cite{Maclaurin18b}, we have demonstrated that such times are typically of exponential length, so there is plenty of time for the Lyapunov exponent to converge to the expected value.
It follows that $\lambda$ corresponds to the long-time average of the right-hand of Eq. (\ref{eq: QSS log phase}). Assuming that the system is ergodic, we can replace the time average by an ensemble average with respect to the Wiener processes. Given
\[\E\left [\sum_{m,n\in \Gamma}{\mathcal F}'_m(\phi)B_{mn}dW_n(t)\right ]=0,\]
for the Ito stochastic process, it follows that
\begin{align}
\lambda_{\rm QSS} &=\epsilon \E\left [\sum_{m,n\in \Gamma}{\mathcal F}'_m(\phi)\tilde{A}_{mn}{\mathcal F}'_n(\phi)\right ]<0,
\end{align} 
provided ${\mathcal F}_n(\phi)$ is not a constant (since the matrix $\tilde{A}$ is negative definite). Since the Lyapunov exponent is then negative definite, we infer that the population of phases evolving according to Eq. (\ref{dphase}) synchronize, in the sense that
\[\lim_{t\rightarrow \infty} [\theta_j(t)-\theta_1(t)]=0,\quad \mbox{for all } j=1,\ldots,M.
\]
Moreover, assuming that a stationary density $P_s(\phi)$ exists, with $P_s(\phi)\approx 1/(2\pi)$ in the weak noise regime, then we can approximate the expectation by an integral around the limit cycle:
\begin{align}
\label{Lyp1}
\lambda_{\rm QSS} &=\epsilon \int_0^{2\pi}\left [\sum_{m,n\in \Gamma}{\mathcal F}'_m(\phi)\tilde{A}_{mn}{\mathcal F}'_n(\phi)\right ] \frac{d\phi}{2\pi}<0.
\end{align} 
This then implies that if we had included $O(\epsilon)$ contributions to the drift, then these would yield a total derivative in the phase difference equation, which would vanish when averaged around the limit cycle. 

In conclusion, we have established that under the diffusion approximation, a population of identical, stochastic hybrid limit cycle oscillators will phase synchronize when driven by a common switching environment in the fast switching limit. This then raises the issue as to whether or not this ensures synchronization of the corresponding population of PDMPs evolving according to the exact dynamics of Eq. (\ref{pdmp}).
In particular, the QSS approximation is only intended to be accurate over timescales that are longer than $O(\epsilon)$. Hence, it is not clear to what extent the above Lyapunov exponent $\lambda_{\rm QSS}$ is accurate, because it is obtained from averaging the noise fluctuations over infinitesimally small timescales. Indeed, it follows from the smoothness of the functions $\lbrace F_m \rbrace$ that the infinitesimal of the exact isochronal phase $\theta_i=\Theta(x_i)$ satisfies $d\theta_i\sim O(dt)$, whereas under the diffusion approximation $d\theta_i\sim \sqrt{\epsilon}O(dt^{1/2})$. This implies that no matter how small we take $\epsilon$, the isochronal phase will never exhibit the relatively large fluctuations over very small timescales that is characteristic of  SDEs.

 In fact the above issue also raises some questions about the conventional approach to phase synchronization in stochastic differential equations. We are not aware of any application of SDEs that is intended to be accurate over infinitely short timescales (i.e. for infinitely high frequencies): in practice, there is always a very short timescale over which the noise is not white, but highly correlated. For example, in stochastic models of stock price fluctuations, this timescale must be at least as long as the time it takes for for the central computer to process a single trade. However, the Lyapunov coefficient that one obtains from the conventional stochastic phase synchronization analysis derives from averaging over these infinitesimally fine fluctuations. In fact, one finds that phase synchronization still occurs in the case of more realistic forms of environmental noise \cite{Goldobin10}. 

\setcounter{equation}{0}

\section{Stochastic hybrid phase equation}

In this section, we consider an alternative approach to deriving the Lyapunov exponent, which avoids the need for the QSS diffusion approximation. The method involves analyzing the PDMP for the exact isochronal phase defined according to
\begin{equation}
\theta_i(t)=\Theta(x_i(t)),
\end{equation}
where $x_i(t)$ now evolves according to the exact PDMP (\ref{pdmp}), rather than the approximate SDE (\ref{dep2}). 

\subsection{Exact PDMP for the isochronal phase}

Suppose that there is a finite sequence of jump times $\{t_1,\ldots t_r\}$ within the time interval $(0,T)$ and let $n_k$ be the corresponding discrete state in the interval $(t_k,t_{k+1})$ with $t_0=0$, see Section II. Introducing the set
\[{\mathcal T}=[0,T]\backslash \cup_{k=1}^r \{t_k\},\]
it follows that Eq. (\ref{pdmp}) holds for all $t\in {\mathcal T}$.
Hence, using the chain rule for $t\in {\mathcal T}$.
\begin{align}
\frac{d\theta_i}{dt}&=\nabla \Theta(x_i)\cdot F_n(x_i)\nonumber \\
&=\overline{\omega}+ \nabla \Theta(x_i)\cdot[F_n(x_i)-\overline{F}(x_i)]\nonumber \\
&=\overline{\omega}+{\mathcal H}_n(x_i),
\label{SHSp}
\end{align}
where
\begin{equation}
{\mathcal H}_n(x):= \nabla \Theta(x)\cdot[F_n(x)-\overline{F}(x)].
\end{equation}
We will use Eq. (\ref{SHSp}) to derive a more accurate expression for the Lyapunov exponent that has the same sign as Eq. (\ref{Lyp1}) but differs in explicit form. The direct method has a number of other advantages. First, it is more intuitive to preserve the piecewise nature of the stochastic dynamics, rather than replacing it by a continuous Markov process. Indeed, it is not clear which aspect of the PDMP the Brownian motion $W_n(t)$ corresponds to. (In appendix B, we explicitly identify a random variable that corresponds to the Brownian motion in the QSS SDE. This then allows us to identify over what timescale the QSS approximation is accurate.) Second, since Eq. (\ref{SHSp}) is exact, it is possible to numerically  solve for $\theta_i(t)$ outside the fast switching regime, and thus determine how phase synchronization varies as $\epsilon$ is increased.

For sufficiently small $\epsilon$, there is a high probability that the environmental state switches multiple times during one period $\Delta_0$. Hence, although ${\mathcal H}_n(x_i)$ is not necessarily $O(\epsilon)$, it only applies for a small time interval before switching, and the accumulative effect of the perturbation over one cycle remains small. This suggests that we can set $x_i=\Phi(\theta_i)$ on the right-hand side of Eq. (\ref{SHSp}), which yields the closed PDMP for the phase:
 \begin{equation}
\label{pdmpq}
\frac{d\theta_i}{dt}=\overline{\omega}+{\mathcal F}_n(\theta_i),\quad t\in {\mathcal T},\quad N(t)=n,
\end{equation}
with ${\mathcal F}_n(\theta)={\mathcal H}_n(\Phi(\theta))$. The corresponding probability density $p_n({\bm \theta},t) $ evolves according to the CK equation 
\begin{equation}
\label{CKphase}
\frac{\partial p_n}{\partial t}=-\sum_{i=1}^M\frac{\partial }{\partial \theta_i} \left [(\overline{\omega}+{\mathcal F}_n(\theta_i))p_n({\bm \theta},t)\right ] +\frac{1}{\epsilon}\sum_{m\in \Gamma}A_{nm}p_m({\bm \theta},t),
\end{equation}
One way to proceed, by analogy with the analysis of SDEs \cite{Nakao07}, would be to consider a pair of oscillators, introduce slow phase variables $\varphi_j=\theta_i-\omega_0t$ and average over a single period of the limit cycle. One could then attempt to find the steady-state probability density for the resulting phase difference, and establish that the phase difference is localized around zero. However, it is difficult to make this approach rigorous, and finding the stationary solution of the CK equation for a PDMP is non-trivial.  

\subsection{Lyapunov exponent}

Therefore, we will proceed by considering a pair of isochronal phases $\theta_1(t)$ and $\theta_2(t)$ evolving according to Eq. (\ref{SHSp}). Set $\theta_1=\phi$, $\theta_2=\phi+\psi$ and define $y=\log \psi$ as before. (Without loss of generality we take $\theta_2(0) > \theta_1(0)$.) Exploiting the fact that the inter-switching times $\Delta t_k=t_{k+1}-t_k$ are exponentially distributed with an $O(1/\epsilon)$ rate, we Taylor expand $\Delta y_k :=y(t_{k+1})-y(t_k)$ to second order in $\Delta t_k$.
 \begin{align*}
\Delta y_k = \psi(t_k)^{-1} \Delta \psi_{k} - \frac{1}{2}\psi(t_k)^{-2} \Delta \psi_{k}^2+ O(\Delta\psi_k^3),
\end{align*}
where $\Delta \psi_k=\psi(t_{k+1})-\psi(t_k)$, and
\begin{eqnarray*}
\Delta \psi_k &=&\Delta t_k  [{\mathcal H}_{n_k}(x_2)-{\mathcal H}_{n_k}(x_1)] \\
&&\quad +\frac{1}{2}\Delta t_k^2 [{\mathcal Q}_{n_k}(x_2)-{\mathcal Q}_{n_k}(x_1)]+O(\Delta t_k^3),\nonumber
\end{eqnarray*}
with $x_j$ evaluated at time $t=t_k$, 
\begin{equation}
{\mathcal Q}_n(x) =   \nabla\Theta(x) \cdot \big\lbrace J_n(x)F_n(x)\big\rbrace+ \big\lbrace J_\Theta(x)F_n(x)\big\rbrace \cdot F_n(x),
\end{equation}
and $J_{\Theta}$ is the Hessian of the isochronal phase map, and $J_n$ is the Jacobian of $F_n$.
Hence, we have
 \begin{eqnarray}
\Delta y_k &=& \psi(t_k)^{-1}  [{\mathcal H}_{n_k}(x_2)-{\mathcal H}_{n_k}(x_1)] \Delta t_k\label{Delta yk decomposition} \\
&&- \frac{1}{2}\psi(t_k)^{-2}  [{\mathcal H}_{n_k}(x_2)-{\mathcal H}_{n_k}(x_1)] ^2\Delta t_k^2\nonumber \\
&& +\frac{1}{2}\psi(t_k)^{-1} [{\mathcal Q}_{n_k}(x_2)-{\mathcal Q}_{n_k}(x_1)]\Delta t_k^2+ O(\Delta t_k^3).\nonumber
\end{eqnarray}

Our goal is to understand the asymptotics of the rate of increase of $y$ with respect to time, i.e. the typical value of 
\begin{equation*}
\frac{\sum_{k=0}^K \Delta y_k}{\sum_{k=0}^K \Delta t_k},
\end{equation*}
for large $K$. We will show that the second term in the decomposition \eqref{Delta yk decomposition} dominates the numerator of the above fraction. It is not immediately obvious why this should be the case, particularly since the second term is itself asymptotically small in $\epsilon$, since $\Delta t_k^2$ scales as $\epsilon^2$. In fact over short time scales, the fluctuations due to $\psi(t_k)^{-1}  [{\mathcal H}_{n_k}(x_1)-{\mathcal H}_{n_k}(x_2)] \Delta t_k$ are dominant. However, we will see that the reason that they are not dominant in the long time average is that their mean is zero, and the fluctuations decorrelate exponentially quickly in time. 

The first step is to note that for sufficiently small $\epsilon$ the system will switch many times during one period $\overline{\Delta}$, while $\psi$ will hardly change. To this end, introduce the cycle times $\tau_p=p\overline{\Delta}$ and set $\hat{\Delta} y_p =y(\tau_{p+1})-y(\tau_p)$. If there are an average of ${\mathcal N}$ jumps during one cycle, we then have
\begin{equation}
\label{dpp}
\hat{\Delta} y_p \simeq \sum_{l=k}^{l=k+{\mathcal N} }\Delta y_l.
\end{equation}
Recall that $\Delta t_k$ is exponentially distributed with rate $\lambda_n \epsilon^{-1}$ when the current discrete state is $n$. This means that 
\[\E[ \Delta t_k |n_k =n]=\frac{\epsilon}{\lambda_n},\quad \E[\Delta t_k^2|n_k=n]=\frac{2\epsilon^2}{\lambda_n^2}.\]
Noting that $\E[ \Delta t_k |n_k =n]$ is the mean waiting time in state $n$, over any particular time interval $\delta t$, we can estimate that the number of jumps 
from state $n$ to some other state per cycle to be $\delta t \rho_n \lambda_n/\epsilon $ (obtained by dividing the expected time spent in state $n$ by the average time it takes to leave state $n$). Summing over $n$, the total number of jumps is approximately 
\begin{equation}\label{eq: total jumps}
\delta t \epsilon^{-1}\sum_{n\in \Gamma}\rho_n \lambda_n.
\end{equation}
In particular, over the course of one cycle, we expect that the total number of jumps is approximately 
\begin{equation}
{\mathcal N}\approx \frac{\overline{\Delta}}{\epsilon}\sum_{n\in \Gamma}\rho_n\lambda_n.\label{eq: number jumps}
\end{equation} 
Since $x_1,x_2$ will not be static over the limit cycle, we must further partition the set of $\mathcal{N}$ jumps into blocks of $\M$ jumps, for $1\ll \M \ll \epsilon^{-1}$, such that over these $\M$ jumps $x_1$ and $x_2$ (as well as $\psi$) are approximately constant. (The motivation for this choice of scaling is as follows:  $\M$ must be small enough that $x_1$ and $x_2$ do not change substantially over all $\M$ jumps, but large enough that the number of jumps that the system makes to each state is determined by $\rho$.) It follows from Eq. \eqref{eq: total jumps} that the total elapsed time over $\M$ jumps is approximately
 \begin{equation}\label{eq: tau M}
\tau_{\M}:=\sum_{l=k}^{l=k+\M} \Delta t_l \simeq \frac{ \epsilon \M}{\sum_{n\in \Gamma}\rho_n\lambda_n}.
 \end{equation}
Furthermore, the number of these $M$ jumps that enter state $n$ is approximately $\mathcal{M}\rho_n \lambda_n / \sum_{m\in\Gamma}\rho_m\lambda_m$.

The next step is to evaluate the sum of each of the terms on the right-hand side of Eq. (\ref{Delta yk decomposition}) over $\M$ jumps. Because there are many jumps, the Law of Large Numbers implies that the sum can be approximated by its average. (We discuss this in more detail in the following section.) For the sake of illustration, let us focus on the second term and consider the  following summation:
 \begin{align*}
&\sum_{l=k}^{l=k+\M }\psi(t_l)^{-2}  [{\mathcal H}_{n_l}(x_2)-{\mathcal H}_{n_l}(x_1)] ^2\Delta t_l^2\\
&\simeq \sum_{l=k}^{l=k+\M }\mathbb{E}\bigg[\psi(t_l)^{-2}  [{\mathcal H}_{n_l}(x_2)-{\mathcal H}_{n_l}(x_1)] ^2\Delta t_l^2\bigg]\\
&\simeq  \frac{\M}{\sum_{m\in \Gamma} \lambda_m \rho_m}\psi(t_k)^{-2} \sum_{m\in \Gamma}\rho_m \lambda_m \lbrace{\mathcal H}_{m}(x_2)-{\mathcal H}_{m}(x_1)\rbrace^2\\ &\qquad \times\mathbb{E}\big[\Delta t_l^2 | n_l = m\big]\\
&=  2\epsilon^{2}\frac{\M}{\sum_{m\in \Gamma} \lambda_m \rho_m}\psi(t_k)^{-2} \sum_{m\in \Gamma}\frac{\rho_m}{ \lambda_m}  \lbrace{\mathcal H}_{m}(x_2)-{\mathcal H}_{m}(x_1)\rbrace^2\\
&\simeq 2\epsilon \tau_{\M}\psi(t_k)^{-2} \sum_{m\in \Gamma}\frac{\rho_m}{ \lambda_m} [{\mathcal H}_{m}(x_2)-{\mathcal H}_{m}(x_1)] ^2.
 \end{align*}
Ignoring transverse amplitude fluctuations, which can be justified using methods from \cite{Maclaurin18b}, we have
\[x_1\simeq \Phi(\phi),\quad x_2\simeq \Phi(\phi+\psi)\simeq x_1+\Phi'(\phi)\psi.\]
Therefore,
\[{\mathcal H}_{m}(x_2)-{\mathcal H}_{m}(x_1)\simeq  {\mathcal H}_{m}'(\Phi(\phi))\Phi'(\phi)\psi.\]
and
\[
\psi(t_k)^{-1} \lbrace{\mathcal H}_{m}(x_2)-{\mathcal H}_{m}(x_1)\rbrace \simeq \mathcal{F}_m'(\phi_k),
\]
where $\mathcal{F}_m(\phi) = \mathcal{H}_m\big(\Phi(\phi)\big)$.
Hence, 
 \begin{align}
\tau_{\M}^{-1}\sum_{l=k}^{l=k+\M }\psi(t_l)^{-2}  [{\mathcal H}_{n_l}(x_1)-{\mathcal H}_{n_l}(x_2)] ^2\Delta t_l^2\nonumber \\
 \simeq 2\epsilon \sum_{m\in \Gamma}\frac{\rho_m}{ \lambda_m}\mathcal{F}'_{m}\big(\phi_k\big)^2
 \label{eq: first derivative}
 \end{align}
We similarly find that
\begin{align}
\tau_{\M}^{-1} \sum_{l=k}^{l=k+\M}\psi(t_l)^{-1} [{\mathcal Q}_{n_l}(x_1)-{\mathcal Q}_{n_l}(x_2)]\Delta t_l^2 \nonumber \\
\simeq 2\epsilon \sum_{m\in \Gamma}\frac{\rho_m}{ \lambda_m}\frac{d}{d\phi_k}\mathcal{Q}_{m}\big(\Phi(\phi_k)\big).
\label{eq: second derivative}
\end{align}

If we apply the same analysis to the summation of the linear term in $\Delta t_k$ in Eq. (\ref{Delta yk decomposition}), we notice that its expectation is approximately zero. That is,
 \begin{align*}
&\mathbb{E}\bigg[ \sum_{l=k}^{l=k+\M }\psi(t_l)^{-1}  [{\mathcal H}_{n_l}(x_1)-{\mathcal H}_{n_l}(x_2)] \Delta t_l\bigg]\\
&\simeq  \frac{\M}{\sum_{m\in \Gamma} \lambda_m \rho_m}\psi(t_k)^{-1} \sum_{m\in \Gamma}\rho_m \lambda_m \lbrace{\mathcal H}_{m}(x_1)-{\mathcal H}_{m}(x_2)\rbrace\\ &\quad \times\mathbb{E}\big[\Delta t_l | n_l = m\big]\\
&=  \epsilon\frac{\M}{\sum_{m\in \Gamma} \lambda_m \rho_m}\psi(t_k)^{-1} \sum_{m\in \Gamma}\rho_m  \lbrace{\mathcal H}_{m}(x_1)-{\mathcal H}_{m}(x_2)\rbrace \\
&=0,
 \end{align*}
since $\sum_{m\in \Gamma}\rho_m \mathcal{H}_m(x_i)= 0$. This suggests that its contribution to the Lyapunov exponent is negligible over long times. A more rigorous analysis indeed establishes that (see Appendix B)
\begin{equation}\label{eq: linear terms in delta t}
\frac{1}{ q\bar{\Delta}}\sum_{l=0}^{q\mathcal{N}}\psi(t_l)^{-1}  [{\mathcal H}_{n_l}(x_1)-{\mathcal H}_{n_l}(x_2)] \Delta t_l \to 0
\end{equation}
as $q\to\infty$, as long as both oscillators stay close to the limit cycle.

Finally, we can take the rate of change of the phase to be approximately that of the deterministic system, over one course of the limit cycle. That is, we assume that for all $t\in [\tau_p, \tau_{p+1}]$, $\phi(t) - \phi(\tau_p) \simeq \bar{\omega}(t-\tau_p)$. Since, by assumption, $\mathcal{M} \ll \epsilon^{-1}$, the expressions in \eqref{eq: first derivative} and \eqref{eq: second derivative} are excellent approximations to the derivative with respect to time. Combining the above results, and choosing $k$ to be such that $\tau_p \simeq t_k$,
\begin{align*}
&y(\tau_{p+1})-y(\tau_p) = \int_{\tau_p}^{\tau_{p+1}} y'(t) dt \\
&\simeq \epsilon \int_{\tau_p}^{\tau_{p+1}}\bigg\lbrace \sum_{m\in \Gamma}\frac{\rho_m}{ \lambda_m}\bigg( \frac{d}{d\phi}\mathcal{Q}_m\big(\Phi(\phi_t)\big) -\mathcal{F}'_{m}(\phi_t)^2\bigg) \bigg\rbrace dt \\
&\simeq \frac{\epsilon \bar{\Delta} }{2\pi}\int_0^{2\pi} \bigg\lbrace \sum_{m\in \Gamma}\frac{\rho_m}{ \lambda_m}\bigg( \frac{d}{d\alpha}\mathcal{Q}_m\big(\Phi(\alpha)\big) -\mathcal{F}'_{m}(\alpha)^2\bigg) \bigg\rbrace d\alpha,
\end{align*}
upon a change of variable. Now
\[
 \int_0^{2\pi}  \frac{d}{d\alpha}\mathcal{Q}_m\big(\Phi(\alpha)\big) d\alpha = \mathcal{Q}_m\big(\Phi(2\pi)\big) -  \mathcal{Q}_m\big(\Phi(0)\big) = 0.
\]
due to the periodicity of $\Phi$. We thus find that, taking $q$ revolutions around the limit cycle (which take a time $ q  \bar{\Delta}$)
\begin{multline*}
\frac{1}{q \bar{\Delta}}\big[y(\tau_{q})-y(\tau_1)\big] 
\simeq 
-\frac{\epsilon}{2\pi } \int_0^{2\pi}  \sum_{m\in \Gamma}\frac{\rho_m}{ \lambda_m}\mathcal{F}'_{m}(\alpha)^2 d\alpha.
\end{multline*}
We thus obtain the Lyapunov exponent
\begin{align}
\label{Lyp2}
\lambda &=-\epsilon \int_0^{2\pi} \sum_{n\in \Gamma}\rho_n\lambda_n^{-1}{\mathcal F}'_n(\phi)^2\frac{d\phi}{2\pi}<0.
\end{align}

\subsection{Remarks}
 
 The above analysis establishes that the Lyapunov exponent $\lambda$ obtained from the exact isochronal phase equation differs significantly from the Lyapunov exponent $\lambda_{\rm QSS}$ obtained under the diffusion approximation, \james{ i.e. }$\lambda \neq \lambda_{\rm QSS}$, with $\lambda_{\rm QSS}$ given by Eq. (\ref{Lyp1}). Since they are both negative definite (assuming ${\mathcal F}_n'(\phi)\neq 0$), they both predict that phase synchronization will occur, but at different rates. The origin of the discrepancy is that $\lambda_{\rm QSS}$ is
 obtained by averaging with respect to noise fluctuations over infinitesimally small timescales that do not occur in the exact PDMP. As we show in appendix B, one can write
 \begin{equation}
 \label{eq: QSS 2}
\lambda_{\rm QSS} \simeq - \lim_{K\to\infty}\frac{1}{2t_K}\mathbb{E}\bigg[\bigg( \sum_{l=0}^{K} \mathcal{F}'_l \Delta t_l\bigg)^2 \bigg],
 \end{equation}
whereas
  \begin{equation}
\lambda \simeq - \lim_{K\to\infty}\frac{1}{2t_K}\mathbb{E}\bigg[ \sum_{l=0}^{K} \big(\mathcal{F}'_l \Delta t_l\big)^2 \bigg].
 \end{equation}
 (Again, the limit should only really be taken up until the time that the system leaves a neighborhood of the limit cycle.)
Our derivation of $\lambda$ is also useful in helping us understand how the QSS approximation works. It is not immediately obvious where the Brownian motion of Eq. (\ref{dep2}) comes from. In appendix B we demonstrate that in Eq. (\ref{Delta yk decomposition}), the term linear in $\Delta t_l$ corresponds to the stochastic integral of the QSS reduction in \eqref{eq: QSS log phase}. More precisely, the probability law of
$ \sqrt{2\epsilon}\int_0^{t_K}\sum_{m,n\in \Gamma}{\mathcal F}'_m(\phi)B_{mn}dW_n(t)$
is very close to the law of
$\sum_{l=0}^K \mathcal{F}'_{n_l}\Delta t_l$.
Indeed, it can be shown that their first two moments are equal to leading order in $\epsilon$. (One could straightforwardly extend the analysis of appendix B to demonstrate that their higher order moments must also converge as $\epsilon \to 0$.) However, the probability laws of the above two random variables are only convergent over timescales much larger than $O(\epsilon)$ and this essentially accounts for the discrepancy in the Lyapunov exponents.

\setcounter{equation}{0}

\section{Example: Radial isochron clock}
In order to illustrate the above general theory, we will consider a particularly simple model of an oscillator based on the complex amplitude equation that arises for a limit cycle oscillator close to a Hopf bifurcation: 
\begin{equation}
\frac{dA}{dt}=(\mu+i\eta)A-(1+i\alpha)|A|^2A ,\quad A \in {\mathbb C}.
\label{ca}
\end{equation}
In polar coordinates $A=r\e^{i\phi}$,
\begin{eqnarray}
\frac{dr}{dt}=r(\mu-r^2),\quad 
\frac{d\phi}{dt}=\eta-\alpha r^2 . 
\end{eqnarray}
This system is also known as a modified radial isochron clock model.
The solution for arbitrary initial data $r(0)=r_0$, $\theta(0)=\theta_0$ is
\begin{subequations}
\label{ric0}
\begin{eqnarray}
r(t)&=&\sqrt{\mu}\left [1+\frac{\mu-r_0^2}{r_0^2}\e^{-2\mu t}\right ]^{-1/2}, \\
\phi(t)&=&\phi_0+\omega t-\frac{\alpha}{2}\log(r_0^2+(\mu-r_0^2)\e^{-2\mu t}),
\end{eqnarray}
\end{subequations}
where $\omega:=\eta-\alpha \mu$ is the natural frequency of the stable limit cycle at $r_0^2={\mu}$. In Cartesian coordinates
\begin{subequations}
\label{ric}
\begin{eqnarray}
\frac{dx}{dt}&:=&F_1(x,y)=\mu x-\eta y-(x^2+y^2)(x-\alpha y)\nonumber  \\ \\
\frac{dy}{dt}&:=&F_2(x,y)=\mu y+\eta x-(x^2+y^2)(y+\alpha x). \nonumber \\
\end{eqnarray}
\end{subequations}

One of the useful features of the radial isochron clock model is that the isochronal phase can be calculated explicitly. Strobing the explicit solution Eq. (\ref{ric0}) at
times $t=n\Delta_0$, we see that
\begin{equation*}
\lim_{n\rightarrow \infty}\phi(n\Delta_0)=\phi_0 -\alpha \ln r_0 \quad \mbox{mod } 2\pi.
\end{equation*}
Hence, we can define an isochronal phase on the whole plane according to
\begin{equation}
\Theta(r,\phi)=\phi-\alpha \ln r .
\label{iso}
\end{equation} 
It follows that the isochrones are logarithmic spirals with $\phi-\alpha \ln r = {\rm constant}$. Now rewrite the phase (\ref{iso}) in Cartesian coordinates,
\begin{equation*}
\Theta(\x)=\tan^{-1}\frac{y}{x}-\frac{\alpha}{2}\log(x^2+y^2),
\end{equation*}
so that
\begin{equation*}
\frac{\partial \Theta}{\partial x}=-\frac{y}{x^2+y^2}-\alpha \frac{x}{x^2+y^2},\quad \frac{\partial \Theta}{\partial y} =\frac{x}{x^2+y^2}-\alpha \frac{y}{x^2+y^2}.
\end{equation*}
On the limit cycle $\Phi(\theta) =\sqrt{\mu}(\cos \theta, \sin \theta)$, so that the components of the PRC are
\begin{subequations}
\begin{align*}
R_x(\theta)&=\frac{\partial \Theta(\Phi(\x))}{\partial x}=\frac{1}{\sqrt{\mu}}\left [-\sin \theta-\alpha \cos \theta\right ] ,\\ R_y(\theta)&=\frac{\partial \Theta(\Phi(\x))}{\partial x}=\frac{1}{\sqrt{\mu}}\left [\cos \theta-\alpha \sin \theta\right ].
\end{align*}
\end{subequations}

\begin{figure*}[t!]
\begin{center}
\includegraphics[width=16cm]{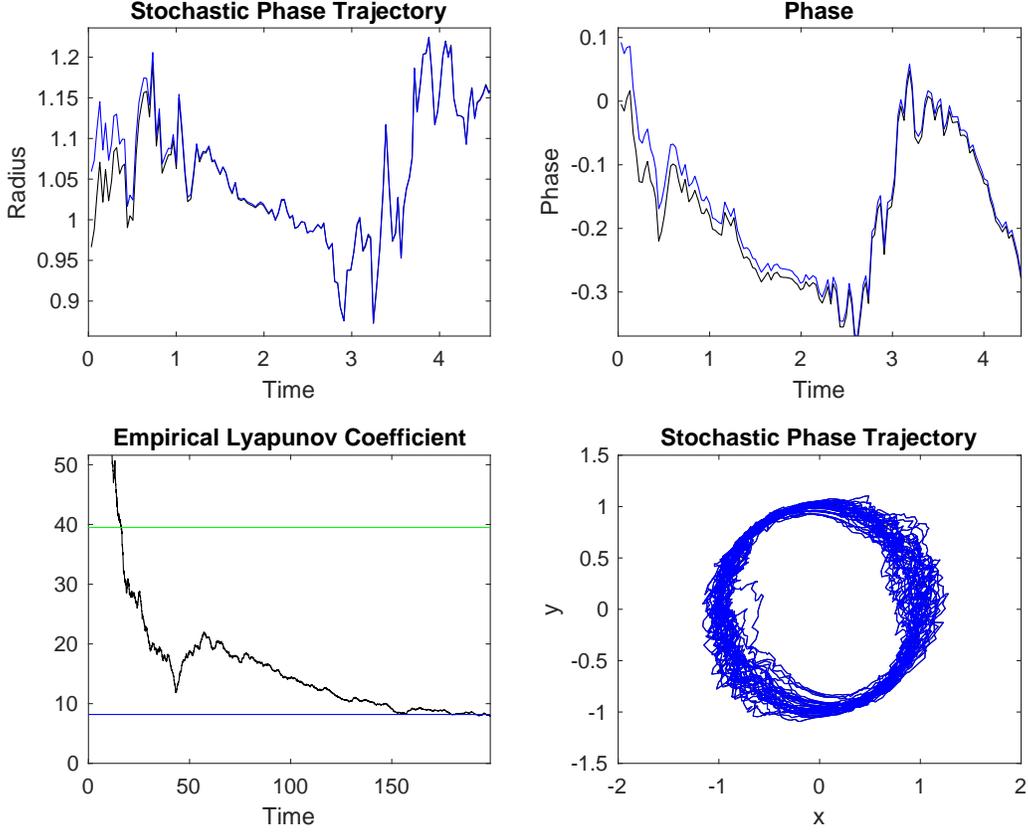}
\caption{\small Synchronization of two radial isochron oscillators with common switching noise. We set $\epsilon = 0.01$. The radii of the initial conditions differ by $0.1$, and the phases of the initial condition differ by $0.1$. In the top left we plot the radii of the two oscillators, and in the top right we plot the angle (one oscillator is plotted in black, and one in blue). In the bottom right we plot the $x-y$ trajectories of the oscillators: one in black, and one in blue. In the bottom left, we plot $-\epsilon^{-1}t^{-1}\log\big| \phi^1(t) - \phi^2(t) \big|$ in black. We plot $\epsilon^{-1}\lambda_{QSS}$ in green, and $\epsilon^{-1}\lambda$ in blue. It can be seen that $\lambda$ is the correct Lyapunov exponent.}
\label{isochrone}
\end{center}
\end{figure*}

Given the deterministic model, the next step is to specify the corresponding PDMP for a single oscillator. One possibility is to assume that one or more of the coefficients switch. For example, in polar coordinates we could take
\begin{eqnarray}
\frac{dr}{dt}=r(\mu_n-r_i^2),\quad 
\frac{d\phi}{dt}=\eta_n-\alpha r^2 . 
\label{ca1}
\end{eqnarray}
That is, both the amplitude and angular frequency of each oscillator switch.
In Cartesian coordinates we have
\begin{subequations}
\label{ric2}
\begin{eqnarray}
\frac{dx}{dt}&:=&F_{n,1}(\x)=\mu_n x-\eta_n y-(x^2+y^2)(x-\alpha y)\nonumber  \\ \\
\frac{dy}{dt}&:=&F_{n,2}(\x)=\mu_n y+\eta_n x-(x^2+y^2)(y+\alpha x), \nonumber \\
\end{eqnarray}
\end{subequations}
It immediately follows that the averaged system is given by
\begin{subequations}
\label{ricav}
\begin{eqnarray}
\frac{dx}{dt}&=&\overline{F}_{1}(\x):=\overline{\mu} x-\overline{\eta} y-(x^2+y^2)(x-\alpha y)\nonumber  \\ \\
\frac{dy}{dt}&=&\overline{F}_{2}(\x):=\overline{\mu} y+\overline{\eta} x-(x^2+y^2)(y+\alpha x), \nonumber \\
\end{eqnarray}
\end{subequations}
where
\begin{equation}
\overline{\mu}=\sum_{n \in \Gamma}\mu_n \rho_n,\quad \overline{\eta}=\sum_{n \in \Gamma}\eta_n \rho_n
\end{equation}
The corresponding natural frequency is $\overline{\omega}=\overline{\eta}-\alpha \overline{\mu}$. Also note from Eq. (\ref{G}) that
\begin{subequations}
\begin{align}
G_{n,1}(\x)=(\mu_n-\overline{\mu})x-(\eta_n-\overline{\eta})y,\\
G_{n,2}(\x)=(\mu_n-\overline{\mu})y+(\eta_n-\overline{\eta})x.
\end{align}
\end{subequations}
The phase PDMP (\ref{pdmpq}) for the radial isochron clock takes the form of a simple velocity jump process
\begin{align}
\label{pdmprci}
\frac{d\theta}{dt}&=\omega_n,
\end{align}
since
\begin{align}
{\mathcal F}_n(\theta)&=\sum_{a=1,2}R_a(\theta)\left [F_{n,a}(\Phi(\theta))-\overline{F}_a(\Phi(\theta))\right ]\nonumber \\
&=\left [-\sin \theta-\alpha \cos \theta\right ] \left [(\mu_n-\overline{\mu})\cos \theta-(\eta_n-\overline{\eta})\sin \theta \right ]
\nonumber\\
&+\left [\cos \theta-\alpha \sin \theta\right ]\left [(\mu_n-\overline{\mu})\sin\theta+(\eta_n-\overline{\eta})\cos \theta\right]\nonumber\\
&= (\eta_n-\overline{\eta})-\alpha(\mu_n-\overline{\mu})=\omega_n-\overline{\omega}.
\end{align}
In this case ${\mathcal F}_n'=0$ so Eq. (\ref{Lyp2}) implies that the Lyapunov exponents is zero \james{(to leading order in $\epsilon$)}-- noise-induced synchronization does not occur. This means that one cannot establish that synchronization occurs using the analysis of this paper. Its possible that synchronization might still occur, but the Lyapunov coefficient would be $o(\epsilon)$, and it is not easy to ascertain the synchronization in the numerical results.

A second possibility is to assume that the environments drives the $x$ coordinate with a switching input $I_n$ such that $\bar{I}=\sum_n\rho_nI_n=0$. As a concrete example, we take the evolution in Cartesian co-ordinates to be given by
\begin{subequations}
\label{ric2}
\begin{eqnarray}
\frac{dx}{dt}&=&\mu x-\eta y-(x^2+y^2)(x-\alpha y) +\frac{x}{\sqrt{x^2+y^2}}v_n^1\nonumber  \\ \\
\frac{dy}{dt}&=&\mu y+\eta x-(x^2+y^2)(y+\alpha x) +\frac{y}{\sqrt{x^2+y^2}}v_n^2 \nonumber \\
\end{eqnarray}
\end{subequations}
Here $v_n = (v^1_k,v^2_k)$ is a jump Markov Process, assuming the following $4$ states
\begin{align*}
v_1 &= (2,-1) \; , \; v_2 = (-4,-4)\\ 
v_3 &= (-3,2) \; , \; v_4 = (8.8,7.2),
\end{align*}
with transition matrix
\begin{equation}
\left( 
\begin{array}{c c c c}
0 & 2 & 2.5 & 0.1 \\
1 & 0 & 0.5 & 4 \\
0.5 & 0.7 & 0 & 2 \\
3 & 0.4 & 0.25 & 0
\end{array}
\right).
\end{equation}
The state vector and transition matrix were chosen arbitrarily, except for the normalization
\[
\sum_{k=1}^4 \rho_k v_k = 0,
\]
to ensure that the averaged system supports a deterministic limit cycle. Hence, the averaged system is given by Eqs. (\ref{ric}), whereas the phase PDMP takes the form
\begin{align}
\label{pdmprci2}
\frac{d\theta}{dt}&=\omega -(\sin \theta+\alpha \cos\theta) v_n^1+\left [\cos \theta-\alpha \sin \theta\right ]v_n^2,
\end{align}
In this case ${\mathcal F}_n'\neq 0$ and we expect a pair of hybrid oscillators to synchronize. 

Numerical simulations of a pair of radial isochron oscillators evolving according to Eqs. (\ref{ric2}) with a common environmental drive confirm that synchronization does occur. An example set of results are shown in Fig. \ref{isochrone}. Note, in particular, that the quantity
\[-\frac{1}{ t}\log\big| \phi^1(t) - \phi^2(t) \big|
\] 
converges to the Lyapunov exponent $\lambda$ given by Eq. (\ref{Lyp2}), which was calculated directly from the underlying PDMP, rather than the quasi-steady-state Lyapunov exponent $\lambda_{\rm QSS}$ of Eq. (\ref{Lyp1}).

\section{Discussion}

In this paper we have proved that \james{stable }oscillators subject to a common rapidly-switching noise will synchronize \james{, in the vast majority of cases}. \james{--Since we gave a counterexample in the previous section, we should probably add this small caveat} We have identified the leading order contribution to the Lyapunov exponent, and explained why this is different from the Lyapunov exponent predicted by the quasi-steady-state assumption. These results were shown to be consistent with a simulation of the radial isochron oscillator subject to a common environmental noise.

In more detail, we have seen that the phases of rapidly-switching oscillators converge at a rate of $\exp\big(-\epsilon \lambda t\big)$, where $\epsilon$ is the timescale of the switching. Thus in the limit as the switching gets faster and faster (i.e. $\epsilon \to 0$), the rate of synchronization gets slower and slower. This is not surprising, since we know that in the deterministic $\epsilon = 0$ limit, the oscillators will in general never synchronize if their starting conditions are different. These results are contingent on the two oscillators staying in the attracting neighborhood of the limit cycle. Indeed we have shown in \cite{Maclaurin18b} that the timescale over which the oscillators remain close to the limit cycle scales as $\exp\big( c\epsilon^{-1}\big)$, for a constant $c$. In fact if one were to continue the analysis of this paper, and develop precise error bounds for the probability of the two oscillators synchronizing, then one would find that the smaller $\epsilon$ is, then the more likely it is that the oscillators, once they are almost synchronized, stay almost synchronized. In summary: as $\epsilon \to 0$, the oscillators synchronize at a slower and slower rate, but stay synchronized with a higher and higher probability.

\bigskip
\section*{Acknowledgements}

PCB and JNM were supported
by the National Science Foundation (Grant No. DMS-1613048). 
\bigskip
\renewcommand{\theequation}{A.\arabic{equation}}

\setcounter{equation}{0}

\section*{Appendix A}

The basic steps of the QSS reduction of the population equations (\ref{dep2}) are as follows:
\medskip

\noindent a) Decompose the probability density as
\begin{equation*}
p_n(\x,t)=C(\x,t)\rho_n +\epsilon w_n(\x,t),
\end{equation*}
where $\sum_{n} p_n(\x,t) =C(\x,t)$ and $\sum_{n} w_n(\x,t)=0$. Substituting into Eq. (\ref{CKH}) yields
\begin{eqnarray*}
 \rho_n \frac{\partial C}{\partial t}+\epsilon \frac{\partial w_n }{\partial t}&=&-\sum_{i=1}^M\nabla_i \cdot \left (F_n(x_i)[\rho_n C+\epsilon w_n]\right )\\
 &&\quad +\frac{1}{\epsilon}\sum_{m \in \Gamma} A_{nm} [\rho_mC+\epsilon w_m]\nonumber
\end{eqnarray*}
Summing both sides with respect to $n$ then gives
\begin{equation}
\label{BBC}
\frac{\partial C}{\partial t}=-\sum_{i=1}^M \left \{\nabla_i \cdot \left [\overline{F}(x_i)C \right ]+\epsilon \sum_{n\in \Gamma} \nabla_i \cdot \left [ F_n(x_i)w_n\right ].\right \}
\end{equation}
\medskip

\noindent b) Using the equation for $C$ and the fact that ${\bf A}\rho= 0$, we have
\begin{eqnarray*}
 &&\epsilon \frac{\partial w_n}{\partial t}=\sum_{m \in \Gamma} A_{nm} w_m\\
 &&- \rho_n\sum_{i=1}^M\left \{ \nabla_i \cdot \left  [F_n(x_i) C \right ]- \nabla_i \cdot \left [\overline{F}(x_i)C\right ] \right \}\\
&&-\epsilon \sum_{i=1}^M\left \{\left [\nabla_i \cdot \left (F_n(x_i)\omega_n\right )- \rho_n\sum_{m\in \Gamma} \nabla_i \cdot \left [ F_m(x_i)w_m\right ]
 \right ]\right \}.
\end{eqnarray*}

\noindent c) Introduce the asymptotic expansion
\[w_n\sim {w}_n^{(0)}+\epsilon {w}_n^{(1)}+\epsilon^2 {w}_n^{(2)}+\ldots\]
and collect $O(1)$ terms:
\begin{eqnarray*}
 &&\sum_{m \in \Gamma} A_{nm}  w^{(0)}_m = \rho_n \sum_{i=1}^M \nabla_i \cdot \left  [ F_n(x_i)-\overline{F}(x_i)]C\right ].\nonumber 
\end{eqnarray*}
The Fredholm alternative theorem show that this has a solution, which is unique on imposing the condition $\sum_n w^{(0)}_n(\x,t)=0$:
\begin{eqnarray}
\label{w0}
 &&w^{(0)}_n(\x)=\sum_{m \in \Gamma} A^{\dagger}_{nm} \rho_m \sum_{j=1}^N\left ( \nabla_j \cdot \left  [ F_m(x_j)-\overline{F}(x_j)]C\right ] \right ).\nonumber 
\end{eqnarray}
where ${\bf A}^{\dagger}$ is the pesudo-inverse of the generator ${\bf A}$.
\medskip

\noindent d) Combining Eqs. (\ref{w0}) and (\ref{BBC}) shows that $C$ evolves according to the Fokker-Planck (FP) equation
\begin{eqnarray*}
\frac{\partial C}{\partial t}&=&- \sum_{i=1}^M\nabla_i \cdot \left [\overline{F}(x_i) C\right ]-\epsilon \sum_{i,j=1}^M \sum_{n,m\in \Gamma}A_{nm}^{\dagger} \rho_m  \\
&& \times \nabla_i\cdot      \left ( F_n(x_i)\nabla_j \cdot \left  [ F_m(x_j)-\overline{F}(x_j)]C\right ] \right ).      
\end{eqnarray*}
Using the fact that $\sum_mw_n=0$, this can be rewritten in the Stratonovich form (\ref{zFP0}).
One typically has to determine the pseudo-inverse of ${\bf A}$ numerically.

\renewcommand{\theequation}{B.\arabic{equation}}

\setcounter{equation}{0}

\section*{Appendix B}

In this appendix we show that the $O(\Delta t_k)$ term in Eq. (\ref{Delta yk decomposition}) does not contribute
to the Lyapunov exponent in the long time limit. First, recall that
\[
 \psi(t_l)^{-1}  [{\mathcal H}_{n_l}(x_1)-{\mathcal H}_{n_l}(x_2)] \simeq \mathcal{F}'_{n_l}(\phi).
\]
From now on, we drop the dependence of $\mathcal{F}$ on $\phi$ to simplify notation.
We start by understanding the leading order second moment of $\M$ terms, i.e.
\[ 
V_{\M} := \mathbb{E}\bigg[\bigg( \sum_{l=k}^{l=k+\M } \mathcal{F}'_{n_l} \Delta t_l\bigg)^2 \bigg] .
\]
(An important reason that we do this is that, as noted in Sect. IV.C, the rate of change of this term with respect to time yields the Lyapunov exponent of the quasi-steady-state reduction.)
$\M$ will be taken to be large enough that the system has switched sufficiently many times for the quasi-steady-state approximation to be accurate, but $\M$ is also taken to be sufficiently small that $x_1,x_2$ and $\phi$ are approximately constant. In other words, $1 \ll \M \ll \epsilon^{-1}$. To this end, we expand out the square to obtain
\begin{multline*}
\mathbb{E}\bigg[\bigg( \sum_{l=k}^{l=k+\M } \mathcal{F}'_{n_l} \Delta t_l\bigg)^2 \bigg] \\
=\sum_{l=k}^{l=k+\M }\sum_{j=k}^{j=k+\M }  \mathbb{E}\bigg[\Delta t_l \Delta t_j  \mathcal{F}'_{n_l}  \mathcal{F}'_{n_j} \bigg] \\
\simeq 2\sum_{l=k}^{l=k+\M}\sum_{j=0}^{\M} \mathbb{E}\bigg[\Delta t_l \Delta t_{l+j}  \mathcal{F}'_{n_l}  \mathcal{F}'_{n_{j+l}} \bigg].
\end{multline*}
The reason why the above approximation is very accurate is that the number of extra terms obtained through the reindexing is negligible compared to the total number of terms.

Now we can approximate the above equation by using the fact that the correlations between $\Delta t_l \mathcal{F}'_{n_l}$ and $\Delta t_{l+j}  \mathcal{F}_{n_{l+j}}'$ decay exponentially fast in $j$, thanks to the Perron-Frobenius Theorem. Since, as demonstrated in Sect. IV, the mean of $\Delta t_{l+j}\Delta F_{n_{l+j}}$ is $0$, the contribution of the terms with asymptotically large $j$ will be negligible. We now explain these statements in more detail.

Let $P = (P_{km})$ be the matrix with elements $P_{km} = \lambda_m^{-1}W_{km}$. If one knows that the system was in state $m$, and that it has jumped, then it jumps to state $k$ with probability $P_{km}$. The Perron-Frobenius Theorem implies that
\begin{equation}\label{eq: Perron Frobenius P}
\lim_{q\to\infty} P^q = B,
\end{equation}
where $B$ is the rank-$1$ matrix with the $i^{th}$ element of each column equal to $\frac{1}{\sum_{a\in \Gamma}\rho_a \lambda_a}\lambda_i \rho_i$. Furthermore the convergence is exponentially fast, i.e.
\begin{equation}\label{eq: Pq exp}
|| P^q - B || = O\big(\exp(-q)\big),
\end{equation}
for any matrix norm \cite{saloff1997lectures}. 

Now
\begin{multline*}
2\sum_{l=k}^{l=k+\M}\sum_{j=0}^{\M} \mathbb{E}\bigg[\Delta t_l \Delta t_{l+j} \mathcal{F}'_{n_l}  \mathcal{F}'_{n_{j+l}} \bigg] = \\
2\sum_{l=k}^{l=k+\M}\sum_{p,m\in \Gamma} \mathbb{P}\big( n_l = m\big)\mathbb{E}\big[\Delta t_l \mathcal{F}'_{n_l} | n_l =m\big]\\ \times \sum_{j=0}^\infty\big(P^j\big)_{pm} \mathbb{E}\big[\Delta t_{l+j} \mathcal{F}'_{n_{j+l}} | n_{j+l} = p \big] \\
\simeq 2\epsilon^2 \sum_{l=k}^{l=k+\M}\sum_{p,m\in \Gamma}\frac{\mathcal{F}'_p \mathcal{F}'_{m}}{\lambda_m \lambda_p} \mathbb{P}\big( n_l = m\big)\sum_{j=0}^{\M}\big(P^j\big)_{pm}.
\end{multline*}
since, because $x_1,x_2$ change by a negligible amount over this timescale (as $\mathcal{M} = o(\epsilon^{-1})$), if $n_i = n_j$, then $\mathcal{F}'_{n_i} \simeq \mathcal{F}'_{n_j}$.

If we substitute $P^j$ for its limit $B$ in the above, we obtain
\begin{equation*}
\frac{1}{\sum_{a\in \Gamma}\rho_a \lambda_a} \sum_{j=0}^{\M}\sum_{l=k}^{l=k+\M}\sum_{p,m\in \Gamma}\frac{\mathcal{F}'_m \mathcal{F}'_p}{\lambda_m \lambda_p} \mathbb{P}\big( n_l = m\big)\rho_p \lambda_p=0,
\end{equation*}
since $\sum_{p\in\Gamma}\rho_p\Delta \mathcal{F}'_p = 0$. This is what we expect in light of the above discussion, because the system decorrelates after infinitely many jumps,  and the mean is zero.

In light of \eqref{eq: Pq exp}, the above discussion means that for large $\M$ we can extend the summation to $\infty$, without much of a loss of accuracy
\begin{multline}
2\sum_{l=k}^{l=k+\M}\sum_{j=0}^{\M} \mathbb{E}\bigg[\Delta t_l \Delta t_{l+j} \mathcal{F}'_{n_l} \mathcal{F}'_{n_{j+l}} \bigg] \\
\simeq 2\epsilon^2 \sum_{l=k}^{l=k+\M}\sum_{p,m\in \Gamma}\frac{\mathcal{F}'_p \mathcal{F}'_{m}}{\lambda_m \lambda_p} \mathbb{P}\big( n_l = m\big)\sum_{j=0}^\infty\big(P^j \big)_{pm} \\
\simeq \frac{2\epsilon^2 \M}{\sum_{a\in\Gamma}\rho_a\lambda_a}\sum_{p,m\in \Gamma}\frac{\rho_m \mathcal{F}'_p   \mathcal{F}'_{m}}{ \lambda_p} \sum_{j=0}^\infty\big(P^j \big)_{pm}, \label{eq:Expect temp}
\end{multline}
since, as explained above, the number of jumps to state $m$ is approximately $\M\rho_m \lambda_m / \big( \sum_{a\in\Gamma}\rho_a\lambda_a\big)$. 

Now we wish to understand how the above term relates to the pseudo-inverse of $A$ (which occurred in the QSS SDE of Appendix A). In fact we claim that
\begin{equation}\label{eq:A pseudo inverse}
-\sum_{j=0}^\infty \sum_{p,q \in \Gamma} A_{mp}\lambda_p^{-1} \big(P^j\big)_{pq} \rho_q \mathcal{F}'_q =  \rho_m\mathcal{F}'_m,
\end{equation}
which will allow us to establish the identity in \eqref{eq:A pseudo inverse 2}.

Note that the summation on the left converges because $\sum_{p\in \Gamma}A_{mp}\lambda_p^{-1}B_{pq} = 0$ for every $q,m\in\Gamma$, and therefore the convergence is exponential in $j$ thanks to \eqref{eq: Pq exp}. We expand out the left hand side, substituting $A$ and using a truncated summation, i.e. for some $R \in \mathbb{Z}^+$,
\begin{align*}
&-\sum_{j=0}^R \sum_{p,q \in \Gamma} A_{mp}\lambda_p^{-1} \big(P^j\big)_{pq} \mathcal{F}'_q \rho_q   \\
&=-\sum_{j=0}^R \sum_{p,q \in \Gamma} \big(W_{mp} - \lambda_p \delta(m,p)\big)\lambda_p^{-1} \big(P^j\big)_{pq} \mathcal{F}'_q\rho_q \\
&=\sum_{q\in \Gamma}\bigg\lbrace \sum_{j=0}^R \big(P^j\big)_{mq} \mathcal{F}_q'\rho_q - \sum_{j=1}^{R+1} \big(P^j\big)_{mq} \mathcal{F}'_q\rho_q\bigg\rbrace \\
&= \mathcal{F}'_m\rho_m - \sum_{q \in \Gamma}\big(P^{R+1}\big)_{mq} \mathcal{F}'_q\rho_q.
\end{align*}
Now  $\sum_{q\in\Gamma}B_{mq}\mathcal{F}'_q\rho_q = 0$, because, as noted just below \eqref{eq: Perron Frobenius P}, the columns of $B$ are all the same, and $\sum_{q\in\Gamma}\mathcal{H}_q \rho_q = 0$. Since $P^{R+1} \to B$ as $R\to\infty$, when we take $R\to\infty$,
we obtain \eqref{eq:A pseudo inverse}. 

It follows from \eqref{eq:A pseudo inverse} that
\begin{equation}\label{eq:A pseudo inverse 2}
\sum_{j=0}^\infty \sum_{q \in \Gamma} \lambda_p^{-1} \big(P^j\big)_{pq} \rho_q \mathcal{F}'_q = -A^{\dagger}_{pm} \rho_m\mathcal{F}'_m,
\end{equation}
where we recall the pseudo-inverse $A^{\dagger}$, defined in the previous section on the QSS approximation.

Substituting this identity into \eqref{eq:Expect temp}, we find that
\begin{multline}
\sum_{l=k}^{l=k+\M}\sum_{j=0}^{\M} \mathbb{E}\bigg[\Delta t_l \Delta t_{l+j} \mathcal{F}'_{n_l}  \mathcal{F}'_{j+l} \bigg] \\ \simeq -\frac{\epsilon^2 \M}{\sum_{a\in\Gamma}\rho_a\lambda_a}\sum_{p,q\in\Gamma}\mathcal{A}^{\dagger}_{pq} \mathcal{F}'_p \mathcal{F}'_q \rho_q.\label{eq: second moment}
\end{multline}
This implies that
\begin{equation}
\frac{1}{2(t_{k+\M} - t_k)}\mathcal{V}_M \simeq -\epsilon\sum_{p,q\in\Gamma}\mathcal{A}^{\dagger}_{pq} \mathcal{F}'_p \mathcal{F}'_q\rho_q,
\end{equation}
using the expression for $\tau_{\M} = t_{k+\M} - t_k$ in \eqref{eq: tau M}. Just like the end of Sect. IV.B, we find that
\begin{eqnarray}
&&\frac{1}{2t_K}\mathbb{E}\bigg[\bigg( \sum_{l=0}^{K} \mathcal{F}'_l \Delta t_l\bigg)^2 \bigg] \nonumber\\
&&\qquad \simeq -\frac{\epsilon}{2\pi} \int_0^{2\pi} \sum_{p,q\in\Gamma}\mathcal{A}^{\dagger}_{pq} \mathcal{F}'_p(\theta) \mathcal{F}'_q(\theta) \rho_qd\theta .\nonumber\\ &&
\qquad  =\lambda_{QSS}
\end{eqnarray}
We have thus established the claim in \eqref{eq: QSS 2}.

The final step in the proof is to show that
\begin{equation}
\mathbb{P}\bigg(t_{q\M}^{-1} \sum_{l=0}^{q\M } \mathcal{F}'_{n_l} \Delta t_l \geq O(\epsilon)\bigg) \ll 1.
\end{equation}
(Technically, this probability is conditional on both systems staying in a neighborhood of the limit cycle until the time $t_{q\M}$. We have demonstrated that this occurs for very long times, with very high probability, elsewhere \cite{Maclaurin18b}.) 
To this end, we use Chebyshev's Inequality, noting that
\begin{multline}
\mathbb{P}\bigg(t_{q\M}^{-1} \sum_{l=0}^{q\M } \mathcal{F}'_{n_l} \Delta t_l \geq a\bigg) \\ \leq \frac{1}{(at_{q\M})^2}\mathbb{E}\bigg[\bigg(\sum_{l=0}^{q\M } \mathcal{F}'_{n_l} \Delta t_l\bigg)^2 \bigg].\label{eq: Chebyshev}
\end{multline}
We now bound the rate of growth in time of the expectation on the right. The immediate use of this bound will be to show that the probability is negligible. A secondary use is that it will help us understand how this discrete approximation compares to the QSS SDE.

We must thus understand the second order moment of $q\M$ terms, ie.
\begin{align*}
&\mathbb{E}\bigg[\bigg(\sum_{l=0}^{q\M } \mathcal{F}'_{n_l} \Delta t_l\bigg)^2 \bigg] 
=\sum_{l=0}^{q\M}\sum_{j=0}^{q\M}  \mathbb{E}\bigg[\Delta t_l \Delta t_j  \mathcal{F}'_{n_l} \mathcal{F}'_{n_j} \bigg] \\
&\quad \simeq 2\sum_{l=0}^{q\M}\sum_{j=0}^{q\M} \mathbb{E}\bigg[\Delta t_l \Delta t_{l+j} \mathcal{F}'_{n_l}  \mathcal{F}'_{n_{j+l}} \bigg] \\
&\quad = 2\sum_{r=1}^q\sum_{l=(r-1)\M}^{l=r\M}\sum_{j=0}^{q\M} \mathbb{E}\bigg[\Delta t_l \Delta t_{l+j}  \mathcal{F}'_{n_l}  \mathcal{F}'_{n_{j+l}} \bigg].
\end{align*}
Using the estimate derived in \eqref{eq: second moment}, we find that
\begin{widetext}
\begin{align*}
\mathbb{E}\bigg[\bigg( \sum_{l=0}^{q\M }\mathcal{F}'_{n_l} \Delta t_l\bigg)^2 \bigg] -\mathbb{E}\bigg[\bigg( \sum_{l=0}^{(q-1)\M } \mathcal{F}'_{n_l} \Delta t_l\bigg)^2 \bigg] \simeq -\frac{2\epsilon^2 \M}{\sum_{a\in\Gamma}\rho_a\lambda_a}\sum_{p,u\in\Gamma}{\mathcal{A}}_{pu}^{\dagger} \mathcal{F}'_p \mathcal{F}'_u \rho_u.
\end{align*}
Hence, using the expression for the time for $\M$ jumps in \eqref{eq: tau M}, 
\begin{align*}
\frac{1}{\sum_{l=(q-1)\M+1}^{q\M}\Delta t_l}\left \{\mathbb{E}\bigg[\bigg( \sum_{l=0}^{q\M } \mathcal{F}'_{n_l} \Delta t_l\bigg)^2 \bigg] -\mathbb{E}\bigg[\bigg( \sum_{l=0}^{(q-1)\M } \mathcal{F}'_{n_l} \Delta t_l\bigg)^2 \bigg]\right \}  & \simeq -\frac{\sum_{n\in \Gamma}\rho_n\lambda_n}{ \epsilon \M}\frac{2\epsilon^2 \M}{\sum_{a\in\Gamma}\rho_a\lambda_a}\sum_{p,u\in\Gamma}{\mathcal{A}}_{pu}^{\dagger} \mathcal{F}'_{p}  \mathcal{F}'_u\rho_u\\
&= -2\epsilon\sum_{p,u\in\Gamma}{\mathcal{A}}_{pu}^{\dagger}\mathcal{F}'_p\mathcal{F}'_u\rho_u.
\end{align*}
Since the above expression approximates the derivative with respect to time, after re-integrating we find that
\begin{equation}\label{eq: square expectation}
\mathbb{E}\bigg[\bigg( \sum_{l=0}^{q\M } \mathcal{F}'_l \Delta t_l\bigg)^2 \bigg] \simeq -2\epsilon \int_0^{t_{q\M}}\sum_{p,u\in\Gamma}{\mathcal{A}}^{\dagger}_{pu} \mathcal{F}'_p(\theta_{s}) \mathcal{F}'_u(\theta_{s}) \rho_uds.
\end{equation}
\end{widetext}
It then follows from the above expression and \eqref{eq: Chebyshev} that
\begin{equation}
\mathbb{P}\bigg(t_{q\M}^{-1} \sum_{l=0}^{q\M } \mathcal{F}'_{n_l} \Delta t_l \geq a\bigg) = O\bigg(\frac{\epsilon}{a^2 t_{q\M}}\bigg).
\end{equation}
This clearly goes to zero as $q \to \infty$. We have thus justified why the linear terms in $\Delta t_l$ \james{ in \eqref{Delta yk decomposition}} are $o(\epsilon)$ over long periods of time, \james{ and this implies \eqref{eq: linear terms in delta t}}. This is why the quadratic terms in $\Delta t_l$ dominate the linear ones over long periods of time.

\bibliographystyle{plain}

\end{document}